\documentclass [11pt]{article}
\usepackage{amsfonts}
\usepackage{amsmath}
\usepackage{latexsym}
\usepackage{epsfig}
\usepackage{graphicx}
\newcommand{\ignore}[1]{}

\newcommand{\logpower}[1]{(\log n)^{#1}}
\def\bN{\overline{N}}

\def\KS2{{\cal KS}2}

\newcommand{\brac}[1]{\left(#1\right)}
\newcommand{\bfrac}[2]{\brac{\frac{#1}{#2}}}

\newcommand{\set}[1]{\left\{#1\right\}}

\def\cC{{\cal C}}
\def\cE{{\cal E}}

\def\ep{\varepsilon}

\def\e{\varepsilon}

\def\E{{\bf E}}

\def\a{\alpha}
\def\b{\beta}

\def\f{\phi}
\def\g{\gamma}

\def\l{\lambda}

\def\m{\mu}
\def\n{\nu}

\def\Pr{\mathbb{P}}

\def\whp{{\bf whp}}
\def\Whp{{\bf Whp}}

\newtheorem{theorem}{Theorem}[section]
\newtheorem{lemma}[theorem]{Lemma}

\newtheorem{conj}[theorem]{Conjecture}

\newtheorem{claim}[theorem]{Claim}

\newcommand{\proofstart}{{\bf Proof.\hspace{1em}}}

\newcommand{\proofend}{\hspace*{\fill}\mbox{$\Box$}}

\newcommand{\rdup}[1]{\lceil #1 \rceil }

\date{}

\begin{document}

\title{\vspace{-1cm} The game chromatic number of random graphs}

\author{Tom Bohman\thanks{Department of Mathematics, Carnegie Mellon
University, Pittsburgh PA15213. E-mail: tbohman@math.cmu.edu. 
Research supported in part by NSF grant DMS-0401147.}
\and Alan Frieze\thanks{
Department of Mathematics, Carnegie Mellon
University, Pittsburgh PA15213. E-mail: alan@random.math.cmu.edu.
Research supported in part by NSF grant CCF0502793 .}
\and Benny Sudakov \thanks{
Department of Mathematics, Princeton University, Princeton, NJ 08544, and
Institute for Advanced Study, Princeton. E-mail:
bsudakov@math.princeton.edu.
Research supported in part by NSF CAREER award DMS-0546523, NSF grant
DMS-0355497, USA-Israeli BSF grant, Alfred P. Sloan fellowship, and
the State of New Jersey.
}
}
\maketitle
\begin{abstract}
Given a graph $G$ and an integer $k$, two players take turns coloring the vertices of 
$G$ one by one using $k$ colors so that neighboring vertices get different colors.
The first player wins iff at the end of the game all the vertices of $G$ are colored.
The game chromatic number $\chi_g(G)$ is the minimum $k$ for which the first player
has a winning strategy.  In this paper we analyze the asymptotic behavior of this 
parameter for a random graph $G_{n,p}$. We show that with high probability 
the game chromatic number of
$G_{n,p}$ is at least twice its chromatic number but, up to a multiplicative constant, 
has the same order of magnitude. We also study the game chromatic number of random bipartite graphs.
\end{abstract}

\section{Introduction}
Let $G=(V,E)$ be a graph and let $k$ be a positive integer. Consider the following game in which
two players Maker and Breaker take turns coloring the vertices of $G$ with
$k$ colors. Each move consists of choosing an uncolored vertex of the graph and assigning to it
a color from $\{1, \ldots, k\}$ so that resulting coloring is {\em proper}, 
i.e., adjacent vertices get different colors.
Maker wins if all the vertices of $G$ are eventually colored. Breaker wins if at some point in the 
game the current partial coloring cannot be extended to a complete coloring of $G$, 
i.e., there is an uncolored vertex such that each of the $k$ colors appears
at least once in its neighborhood. We assume that Maker goes first (our results will not
be sensitive to this choice). The {\em game chromatic number}
$\chi_g(G)$ is the least integer $k$ for which Maker has a winning
strategy. 

This parameter is well defined, since it is easy to see that Maker always wins if the number of 
colors is larger than the maximum degree of $G$.  Clearly, $\chi_g(G)$ is at least as large as  
the ordinary chromatic number $\chi(G)$, but it can be considerably more.  For example,
let $G$ be a  complete bipartite graph $K_{n,n}$ minus a perfect
matching \(M\) and consider the following strategy for Breaker.
If Maker colors vertex \(v\) with color \(c\) then Breaker
responds by coloring the vertex \(u\) matched with \(v\) in the matching \(M\)
with the same color \(c\).  Note that now \(c\) cannot be used on any other
vertex in the graph.  Therefore, if the number of colors is less than \( n \),
Breaker wins the game. This shows that there are bipartite graphs with 
arbitrarily large game chromatic number and thus there is no upper bound on
$\chi_g(G)$ as a function of $\chi(G)$.

The game was first considered by Brams about 25 years ago in the
context of coloring planar graphs and was described in Martin
Gardner's column \cite{Ga} in Scientific American in 1981. The game remained
unnoticed by the graph-theoretic community until
Bodlaender \cite{Bod} re-invented it. It
has been studied for various classes of graphs in recent years.
Faigle, Kern, Kierstead and Trotter \cite{FKKT} proved that 
the game chromatic number of a forest is at most $4$, and that there are forests
which require that many colors. The game chromatic number of planar graphs
was studied by Kierstead and Trotter \cite{KTr}, who showed that 
for such graphs the game chromatic number is at most $33$. Moreover they proved that 
any graph embeddable on an orientable surface of genus $q$ has                                
game chromatic number bounded by a function of $q$.    
Several additional results on $\chi_g$ and some related parameters were obtained in  
\cite{BK, DZ, GZ, Z1, Kr, NS, CZ, KTu}. For a recent survey see
Bartnicki, Grytczuk, Kierstead and Zhu \cite{BGKZ}.

In this paper, we study the game chromatic number of the random graph \( G_{n,p} \).
As usual, $G_{n,p}$ stands for the probability space of all labeled graphs on $n$ vertices, 
where every edge appears independently with probability $p=p(n)$. We assume throughout the paper that 
the edge probability $p\le 1 - \eta$, where $\eta>0$ is an arbitrarily small, but fixed, constant.
Define $b=\frac{1}{1-p}$ and note that $\log_b x= \frac{\log x}{\log 
b}=(1+o(1))\frac{\log x}{p}$ for all $x\geq 1$ and $p=o(1)$.
Our first result determines the order of magnitude of the game chromatic number of $G_{n,p}$.

\begin{theorem}\label{th1}\ 
\begin{description}
\item[(a)] There exists $K>0$ such that for $\e>0$ and $p \geq (\log n)^{K\e^{-3}} /n$
we have that \whp
\footnote{A sequence of events $\cE_n$ occurs {\em with high
probability} (\whp) if $\lim_{n\to\infty}\Pr(\cE_n)=1$}
$$\chi_g(G_{n,p})\geq (1-\e)\frac{n}{\log_bnp}\, . $$

\item[(b)] If $\a>2$ is a constant, \( K=\max\{\frac{2\a}{\a-1},\frac{\a}{\a-2}\} \) 
and $p \geq (\log n)^K/n$ then \whp
$$\chi_g(G_{n,p})\leq \a \frac{n}{\log_bnp}\, .$$
\end{description}
\end{theorem}

It is natural to compare our bounds with the asymptotic behavior of the
ordinary chromatic number of random graph. It is known by the results of 
Bollob\'as \cite{Bo} and {\L}uczak \cite{Lu}) that \whp\
$\chi(G_{n,p})=(1+o(1)) \frac{n}{2\log_b np}$. Thus our result
shows that the game chromatic number of $G_{n,p}$
is at least twice its chromatic number, but up to a multiplicative constant
has the same order of magnitude.

As already mentioned, there are graphs whose game chromatic number is much larger then
the ordinary one. Our next theorem provides the existence of a large collection of such graphs.
Let $B_{n,p}$ denote the random bipartite graph with two parts of $n$ vertices 
where each of the $n^2$ possible edges appears randomly and independently with
probability $p$. We obtain the following bounds on the game chromatic number of this graph.

\begin{theorem}\label{th2}\ 
\begin{description}
\item[(a)] If $p \geq 2/n$ then 
$$\chi_g(B_{n,p})\geq \frac{n}{10(\log n)(\log_bnp)}\, .$$

\item[(b)] If $\a>2$ is a constant,  \( K=\max\{\frac{2\a}{\a-1},\frac{\a}{\a-2}\} \) 
and $p \geq (\log n)^K/n$ then \whp
$$\chi_g(B_{n,p})\leq \a \frac{n}{\log_bnp}\, .$$
\end{description}
\end{theorem}

The rest of this paper is organized as follows. The next two sections contain proofs of   
lower and upper bounds in Theorem \ref{th1}.  In Section \ref{sec4} we consider the game chromatic number 
of random bipartite graphs and prove Theorem \ref{th2}. 
The last section of the paper contains some concluding remarks and 
open problems. 

Unless the base is specifically mentioned, $\log$ will refer to natural 
logarithms. We often refer to the following Chernoff-type bounds for the tails of binomial
distributions (see, e.g., \cite{AS} or \cite{JLR}). Let $X=\sum_{i=1}^n X_i$ be a sum of 
independent indicator random variables such that $\Pr(X_i=1)=p_i$ and let 
$p=(p_1+\cdots+p_n)/n$. Then
\begin{eqnarray*}
\Pr(X\leq (1-\e)np)&\leq&e^{-\e^2np/2},\\
\Pr(X\geq (1+\e)np)&\leq&e^{-\e^2np/3},\qquad\e\leq 1,\\
\Pr(X\geq \m np)&\leq&(e/\m)^{\m np}.
\end{eqnarray*}

\section{Lower bound on the game chromatic number of $G_{n,p}$}
Suppose that $p \geq (\log n)^{K\e^{-3}}/n$, where $K$ is a sufficiently large constant, and  
that the number of colors $k$ satisfies
$k\leq (1-\e) \frac{n}{\log_bnp}$.
We begin by defining a series of numbers (which will serve as
cut-offs for Breaker's strategy). Let
$$\ell_1=\log_bn-\log_b\log_bnp-10\log_b\log n$$
and note that \(\ell_1\) has been chosen so that we have
\begin{equation}
\label{l_1}
(1-p)^{\ell_1}= \frac{\log_bnp \logpower{10}}{n}=(1+o(1))\frac{\ell_1 \logpower{10}}{n}.
\end{equation}
Set
$$\ell_2 = \frac{\e \ell_1}{20} \ \ \ \ \ \text{ and } \ \ \ \ \
\ell_3=  \frac{ \e^3\ell_1}{2 \cdot 10^6 }.$$
For $S\subseteq [n]$ let $N(S)$ be the neighbors of $S$ which are not in $S$
and let $\bN(S)=[n]\setminus (S\cup N(S))$.

We are now ready to describe Breaker's strategy.
Fix a color $i$. Whenever Maker uses $i$, Breaker will respond by 
using color $i$ in his 
next move. At the beginning Breaker chooses this vertex arbitrarily;
only when $\ell_1$ vertices are colored $i$ will Breaker
choose carefully the next vertex to color. Let $T$ denote the set of
uncolored vertices in $\bN(C_i)$
at the time when the set of vertices $C_i$ that have been colored with
$i$ satisfies $|C_i|=\ell_1$. At this point Breaker identifies a
maximum size independent subset $I_1$ of $T$. When Breaker next uses
color $i$, he will color the vertex $v\in T$ which has as many
neighbors in $I_1$ as possible. After this, $I_1\gets I_1\setminus
N(v)$. When $|I_1|\leq \ell_3$ we say that Breaker has
completed {\em elimination iteration 1}. After completing elimination 
iteration $j$,  Breaker will start a new iteration by identifying the
largest independent set $I_{j+1}$ in the set of uncolored vertices in the
current $\bN(C_i)$ and continue with
the previous strategy. This continues as long as at the start of a new
iteration the set $\bN(C_i)$ contains an independent set of uncolored vertices 
of size at least $\ell_2$. Once $\bN(C_i)$ does not contain any more 
independent sets of size $\ell_2$, from then on Breaker again 
colors arbitrarily with $i$
when desired.

In order to validate Breaker's strategy we must establish a few facts.
We begin by considering the size of \( \bN(C_i) \) when \( |C_i| = \ell_1 \).
\begin{lemma}\label{lem2}
For every subset $S\subseteq [n]$ of size $|S|=
\ell_1$ \whp $~\ell_1 \logpower{9} \leq |\bN(S)| \leq \ell_1 \logpower{11}$.
\end{lemma}
\proofstart
Fix $S$ with $|S|=\ell_1$. The size of $\bN(S)$ is distributed as
the binomial $B(n-\ell_1,(1-p)^{\ell_1})$. Therefore by 
(\ref{l_1}) the expected size of
$\bN(S)$ is $ (1 + o(1)) \ell_1 \logpower{10}$. Thus, it follows from
the Chernoff bounds that
$$\Pr\Big[\,\exists S:\;|S|=\ell_1,\,|\bN(S)|\notin [\ell_1 \logpower{9} ,\ell_1 \logpower{11}]\,\Big] 
\leq
2\binom{n}{\ell_1}e^{-\Theta(\ell_1 \logpower{10})}=o(1).$$
\proofend

Next we note that if the number of uncolored vertices in \( \bN(C_i) \) is sufficiently
large, then Breaker should be able to choose a vertex that reduces the size of
\( I_j \) by a substantial amount.
\begin{lemma}\label{lem3}
\Whp\ there do not exist $S,A,B\subseteq [n]$ such that
\begin{enumerate}
\item $|S|=\ell_1,\,a=|A|\in [\ell_3,3\ell_1],\,|B|\geq b_1=100\e^{-1}\ell_1\logpower{2}$.
\item $A,B\subseteq \bN(S)$ and $A\cap B=\emptyset$.
\item Every $x\in B$ has fewer than $ap/2$ neighbors in $A$.
\end{enumerate}
\end{lemma}
\proofstart
Applying Lemma~\ref{lem2} to bound the 
size of \( \bN(S) \) and using $\ell_3 p=\Omega(K\log \log n)$ we see that
the probability of this event is at most
\begin{eqnarray*}
&&o(1)+\binom{n}{\ell_1}\sum_{a=\ell_3}^{3\ell_1}\binom{\ell_1
\logpower{11}}{a}\binom{\ell_1 \logpower{11}}{b_1}\Pr(B(a,p)\leq ap/2)^{b_1}\\
&\leq&o(1)+n^{\ell_1}\sum_{a=\ell_3}^{3\ell_1}\logpower{11(a+b_1)}e^{-ab_1p/8}\\
&=&o(1).
\end{eqnarray*} 
\proofend

Now we consider the number of elimination iterations.  Define
\[ \ell_2^\prime = \frac{\e \ell_1}{21}, \]
and note that we have
\( \ell_2^\prime < \ell_2 - \ell_3 \).
Each elimination iteration removes at least \( \ell_2^\prime \) vertices
from the set of vertices that can be colored with color \(i\).  Note further
that these sets are disjoint and that each forms an independent set in our graph.
\begin{lemma}\label{lem1}
\Whp\ there do not exist $S,T_1,T_2,\ldots,T_{a_1}\subseteq [n],\,a_1=2000\e^{-2}$
such that
\begin{enumerate}
\item $S,T_1,T_2,\ldots,T_{a_1}$ are pair-wise disjoint independent sets.
\item $|S|=\ell_1$.
\item $|T_i|= \ell_2^\prime,\,i=1,2,\ldots,a_1$.
\item $N(S)\cap T_i=\emptyset,\,i=1,2,\ldots,a_1$.
\end{enumerate}
\end{lemma}
\proofstart
Let $\cE_1$ be the event that such a collection of sets exists.
Then, using $\ell_1/\ell'_2=\Theta(1)$, 
$\frac{\ell_1}{n}=O\big(\frac{\log n}{np}\big)\ll \big(\log 
n\big)^{-K\e^{-3}/2}$ together with 
 (\ref{l_1}) and  Lemma \ref{lem2}, we have
\begin{eqnarray*}
\Pr(\cE_1)&\leq&o(1)+\binom{n}{\ell_1}
\binom{\ell_1 \logpower{11}}{\ell_2^\prime}^{a_1}(1-p)^{a_1\binom{\ell_2^\prime}{2}}\\
&\leq&o(1)+\bfrac{ne}{\ell_1}^{\ell_1}
\brac{\frac{e\, \ell_1\logpower{11}}
{\ell_2^\prime}(1-p)^{(\ell_2^\prime-1)/2}}^{a_1\ell_2^\prime}\\
&\leq&o(1)+\bfrac{ne}{\ell_1}^{\ell_1}
\left(\logpower{12} \Big(\frac{\ell_1 \logpower{10}}{n}\Big)^{\e/42}\right)
^{a_1\e\ell_1/21}\label{force}\\
&\leq&o(1)+ \bfrac{n}{\ell_1}^{\ell_1} \bfrac{\ell_1}{n}^{a_1\e^2\ell_1/1000} \big(\log 
n\big)^{a_1\e \ell_1}\\
&\leq& o(1)+ \left(\frac{\ell_1}{n} \big(\log n\big)^{2000/\e} \right)^{\ell_1}\\
&=& o(1).
\end{eqnarray*}
\proofend

We will complete the proof by showing that most colors are 
used on roughly \( \ell_1 \) vertices.  We will use the following Lemma to bound the
number of colors that are used on significantly more vertices.  We define a fourth 
cut-off
$$\ell_0 =\ell_1 + \frac{ 12\ell_1}{ \ell_3p } \cdot a_1 .$$
Note that \( a_1 \ell_1/ \ell_3\) is a constant and $(1-p)^{\ell_0}= \Omega \left( (1-p)^{\ell_1} \right)$.
\begin{lemma}\label{lem4}
\Whp\ there do not exist pair-wise disjoint sets
$S_1,S_2,\ldots,S_{b_2},U,\,b_2=\frac{n}{\ell_1\logpower{7}}$ such that
\begin{enumerate}
\item $|S_i|= \ell_0$ for $i=1,2,\ldots,b_2$.
\item $|U|=\rdup{n/\log n}$.
\item $|U\cap \bN(S_i)|\leq \ell_1\logpower{8}$ for $i=1,2,\ldots,b_2$.
\end{enumerate}
\end{lemma}
\proofstart
Let $\cE_2$ be the event that such a collection of sets exists.
For every choice of $S_i$ and $U$ the size $|U\cap \bN(S_i)|$ is distributed as 
the binomial $B(n/\log n,(1-p)^{\ell_0})$ with expectation 
$(n/\log n)(1-p)^{\ell_0}=\Omega\big( n(1-p)^{\ell_1}/\log n\big)=\Omega(\ell_1 \logpower{9})$.
Thus, by the Chernoff bound,
\begin{eqnarray*}
\Pr(\cE_2)&\leq&\binom{n}{n/\log n}\brac{\binom{n}{\ell_0}\Pr\big(B(n/\log n,(1-p)^{\ell_0})\leq 
\ell_1\log^8 n\big)}^{b_2}\\
&\leq&n^{n/\log n}n^{ 2 b_2\ell_1}e^{-\Omega(b_2\ell_1\logpower{9})}\\
&=&o(1).
\end{eqnarray*}
\proofend

We now use these Lemmas to complete the proof, assuming 
that the associated low probability events do not occur.
Assume for the sake of contradiction that the
game reaches the point where only $n/\log n$ vertices remain to be
colored.  Let $U$ be the set of uncolored vertices. Let $C_i$
denote the set of vertices colored $i$ at this point and let
$c_i=|C_i|$ for $i=1,2,\ldots,k$. 
Observe that \whp
\begin{equation}\label{indset}
c_i\leq (2+.01\e)\ell_1,\,i=1,2,\ldots,k
\end{equation}
since the right hand side is an upper bound on the size of an independent
set in $G_{n,p}$.

\begin{claim}
\label{cl:bounder}
Let \(i\) be a color such that \( c_i\geq (1+\e/4)\ell_1 \). 
If \( C_i^\prime \) are the first \( \ell_0 \) vertices 
to be colored with color \(i\) then we have
\[ \left| \bN(C_i^\prime) \cap U \right| \le b_1= \frac{ 100 \ell_1\logpower{2}}{\e} \]
\end{claim}
\proofstart

Assume for the sake of contradiction that 
$| \bN(C_i^\prime) \cap U | > b_1 $.
Let $S_t$ be the set of vertices which are
colored $i$ at time $t$. At all times $t$ such that $ |S_t| < \ell_0 $ the
set of uncolored vertices in $\bN(S_t)$ has size at least
$|\bN(C_i^\prime) \cap U| > b_1$.  
Let \( I_j \) be the independent set that is being 
eliminated at time \(t\).  Since \( |I_j| \) is smaller than the independence
number of \( G_{n,p} \),
Lemma~\ref{lem3} implies that Breaker can choose a vertex
that eliminates at least \( |I_j| p /2 \) vertices from \( I_j \).  
Therefore, each 
elimination iteration involves at most $2 
\cdot |I_j|/(\ell_3p/2)<9\ell_1/(\ell_3p)$ 
uses of color $i$. But Lemma~\ref{lem1} implies that 
there are at most $a_1$ elimination iterations.  Therefore, Breaker will complete
all of the elimination iterations before color \(i\) has been 
used \( \ell_0 \) times.  Note that after the elimination process is completed
one can only color at most $\ell_2=\e\ell_1/20$ vertices by color $i$. Therefore,
$$c_i<\ell_1 + 9a_1\ell_1/(\ell_3p)+\e\ell_1/20<(1+\e/4)\ell_1.$$
This is a contradiction.
\proofend

It follows from Claim~\ref{cl:bounder} and Lemma~\ref{lem4} that there at most
\( b_2 = \frac{n}{\ell_1\logpower{7}} \) colors \(i\) such that \( c_i > 
(1+\e/2)\ell_1 \).  Applying this fact together with (\ref{indset})
and $k \leq (1-\ep)n/\ell_1$ we obtain 

$$n-\frac{n}{\log n}=\sum_{i=1}^kc_i\leq
b_2 (2+.01\e)\ell_1+( k- b_2 )(1+\e/4)\ell_1<(1-\e/2)n.$$
This is a contradiction.

\section{Upper bound on the game chromatic number of $G_{n,p}$}

Let $\alpha$ be any constant greater than 2, \(K> \max\{ \frac{ 2 
\alpha}{ \alpha -1} , \frac{ \alpha}{\a -2} \} \), 
$p>\logpower{K}/n$ and let the number of colors be 
$k =\alpha\frac{n}{\log_b np}$.  We begin with Maker's strategy.
Let $\cC=(C_1,C_2,\ldots,C_k)$ be a collection of pair-wise disjoint
sets. Let $\bigcup\cC$ denote $\bigcup_{i=1}^k C_i$. For a vertex $v$ 
let
$$A(v,\cC)= \set{i\in [k]:\;v\mbox{ is not adjacent to any
      vertex of }C_i}. $$
and set
$$a(v,\cC)=|A(v,\cC)|.$$  
Note that \( A(v, \cC) \) is the set of colors that are available at vertex \(v\) when 
the partial coloring is given by the sets in \( \cC \) and \( v \not\in \bigcup \cC\).
Maker's strategy can now be easily defined. Given the current color
classes $\cC$, Maker chooses an uncolored vertex $v$ with the smallest
value of $a(v,\cC)$ and colors it by any available color.

In order to establish that Maker's strategy succeeds \whp, we consider a sequence
of landmarks in the play of the game.  As the game evolves, we let \(u\) denote the
number of uncolored vertices in the graph.  So, we think of \(u\) as running `backward'
from \( n \) to \(0\). Below we define a sequence of thresholds
\( d_0 \ge d_1 \ge \dots \ge d_{r+1} \) and consider the `times' \( u_i \) which are defined 
to be the last  times (i.e. minimum value of \(u\)) for which Maker
colors a vertex for which there are at least \( d_i \) available colors.

We begin with the first landmark, \( u_0 \).
Let
$$\b= k \cdot (np)^{-1/\alpha}=\alpha\frac{n\,(np)^{-1/\alpha}}{\log_b np}\,, 
\hspace{2cm}
\g=\frac{10n\log n}{\b}$$
and
$$B(\cC)=\set{v:\;a(v,\cC)< \b/2}.$$
We begin by showing that with high probability every coloring of the full vertex set 
has the property that 
there are at most 
$\g$ vertices with less than $\b/2$ available colors.

\begin{lemma}\label{lem5}
\Whp, for all collections $\cC$, 
$$|B(\cC)|\leq \g.$$
\end{lemma}
\proofstart
Fix $\cC$. Then for every $v\notin \bigcup \cC$, the number of colors available at $v$ is the 
sum of independent indicator variables $X_i$, where $X_i=1$ if $v$ has no neighbors in $C_i$.
Then $\Pr(X_i=1)=(1-p)^{|C_i|}$ and since $(1-p)^t$ is a convex function we have
\begin{eqnarray*}
\E(a(v,\cC))&=&\sum_{i=1}^k(1-p)^{|C_i|}\\
&\geq&k(1-p)^{(|C_1|+\cdots+|C_k|)/k}\\
&\geq& k(1-p)^{n/k}=\b.
\end{eqnarray*}
It follows from the Chernoff bound that
$$\Pr(a(v,\cC)\leq \b/2)\leq e^{-\b /8}.$$
Thus,
$$\Pr\big(\exists\, \cC~\mbox{with}~|B(\cC)|>\g\big)\leq k^n\binom{n}{\g}e^{-\b\g/8}=o(1).$$
\proofend

Set $u_0$ to be the last time for which Maker colors a vertex with at least 
$d_0=\b/2$ available colors, i.e.,
$$u_0= \min\set{u:\;a(v,\cC_u)\geq d_0=\b/2, \mbox{ for all }v \not\in \bigcup \cC_u},  $$
where $\cC_u$ denotes the collection of color classes when \(u\) vertices remain uncolored. 
It follows from Lemma \ref{lem5} that \whp\ $u_0 \leq \g$ (we apply the Lemma to the 
final coloring). This implies that at some point where the 
number of uncolored vertices is less than $\g$, every vertex still has at least 
$d_0=\b/2$ available colors. In particular, if
$\b/2>\g$ (this happens, e.g., for constant $p$ and $\alpha>2$)
we see that Maker wins the game since no vertex will ever run out of colors.
On the other hand, the proof that Maker's strategy succeeds also for $p=o(1)$
needs more delicate arguments, which we present next.

Since \( u_0 \) is both defined and understood, we are ready to
define \( u_1, \dots, u_{r+1} \).  These landmarks are defined in terms of
thresholds \( d_1, \dots, d_{r+1} \), where \( d_i \) is a lower bound on the number
of colors available at every uncolored vertex (note that \( d_0 = \b/2 \) was 
set above).  We set
$d_{i+1}=d_i-x_i$ where the \( x_i\)'s will be defined below and
$$u_i=\min\set{u:\;a(v,\cC_u)\geq d_i,\mbox{ for all }v \not\in \bigcup \cC_u}$$
We will choose the \( x_i \)'s so that
\begin{equation}\label{sequence}
u_i\geq \log\log n \ \ \ \mbox{ implies } \ \ \ u_{i+1}\leq \frac{u_i}{10}.
\end{equation}
We define $r$ (and hence establish the end of our series of landmarks) 
by $u_r\geq \log\log n>u_{r+1}$. 
Note that if (\ref{sequence}) holds then we have
$$r\leq \log n.$$
We will also ensure that we have
\begin{equation}\label{sequence2}
\sum_{i=0}^rx_i\leq d_0/2.
\end{equation}
If we can choose the \( x_i\)'s so that (\ref{sequence}) and (\ref{sequence2}) 
are satisfied then Maker will succeed in winning the game. Indeed, when there are 
$u_{r}$ uncolored vertices, there are fewer than $10\log\log n$ vertices to color and the 
lists of available colors at these vertices have size at least 
\begin{equation}\label{size}
d_0-\sum_{i=0}^r x_i \geq d_0/2 > \b/4 > 10\log\log n,
\end{equation}
where the lower bound follows from the facts that
$\b \geq \Omega\big((np)^{\frac{\alpha-1}{\alpha}}/\log np\big)$, $np>\logpower{K}$ and 
$K>\frac{2\alpha}{\alpha-1}$.

The key to our analysis (and the choice of \( x_i\)'s) is the following observation.
Between the point when there are
\( u_{i} \) uncolored vertices and the end of the game 
every vertex in \( [n] \setminus \bigcup \cC_{u_{i+1}} \) must
lose at least \( d_{i}-d_{i+1}=x_{i} \) of its available colors.  Indeed, such a vertex \(v\) 
must have at least \( d_{i} \) available colors when there are \(u_i\) vertices uncolored
but has less then \( d_{i+1} \) available colors when \(v\) itself is colored.  This 
implies that the graph induced on 
\( [n] \setminus \bigcup \cC_i \) has at least \( u_{i+1} x_i \) edges.  Before
we proceed, we need another technical 
Lemma, which bounds the number of edges spanned by subsets of $G_{n,p}$.  For each positive 
integer \(s\) define
$$\f=\f(s)=(5ps+\log n)s. $$
\begin{lemma}\label{lem6}
\Whp\ every subset $S$ of $G_{n,p}$ of size $s$ spans at most
$\f=\f(s)=(5ps+\log n)s$ edges.
\end{lemma}
\proofstart
\begin{eqnarray*}
\Pr\big(\exists\,S~\mbox{with}~e(S)>\f\big)&\leq&\sum_{s=2}^n\binom{n}{s}\binom{\binom{s}{2}}{\f}p^{\f}\\
&\leq&\sum_{s=2}^n\brac{\frac{ne}{s}\bfrac{e}{10}^{\log n}}^s\\
&=&o(1).
\end{eqnarray*}
\proofend

We henceforth assume that the low probability events given in Lemmas~\ref{lem5}~and~\ref{lem6}
do not occur.  It follows
from our key observation that we have
\begin{equation}\label{keys}
x_iu_{i+1}\leq \f(u_i)=(5pu_i+\log n)u_i.
\end{equation}

Thus to achieve \eqref{sequence} it suffices to take
$$x_i\geq 10(5pu_i+\log n).$$
But, since $u_i\leq u_0/10^i$ and $u_0\leq \g$, we can take
$$x_i=\frac{5p\g}{10^{i-1}}+10\log n.$$
Checking \eqref{sequence2} we see that we require
$$\sum_{i=0}^rx_i\leq 60p\g+10\logpower{2} $$
to be less than \( d_0/2 \),
and so we need to verify 
\begin{equation}\label{cond}
\frac{600np\log n}{\b}+10\logpower{2}\leq \frac{\b}{4}.
\end{equation}
Note that (\ref{size}) follows immediately from (\ref{cond}).  
Since $np>\logpower{K}$, $K>\max\{\frac{2\a}{\a-1},\frac{\a}{\a-2}\}$ and 
$\alpha>2$ we have that
$$\b \geq \Omega\left(\frac{(np)^{\frac{\alpha-1}{\alpha}}}{\log np}\right)\gg
\max\big\{\logpower{2}, (np)^{1/\alpha} \log n\big\},$$
which implies (\ref{cond}).

\section{Proof of Theorem \ref{th2}}\label{sec4}

We start by proving part (a) of Theorem \ref{th2}.
Suppose that $2/n \leq p\leq 1-\eta$ and the number of colors is at most 
$k =\frac{n}{10(\log n)(\log_b np)}$. Also recall that if 
$p=o(1)$ then $\log_b x=(1+o(1))\frac{\log x}{p}$.

Breaker employs the following strategy.  
He chooses one part of the bipartite graph which we denote by $W_B$ to be {\em Breaker's} side 
and thinks of the opposite part $W_M$ as {\em Maker's} side.  
Loosely speaking, Breaker tries to eliminate the 
coloring possibilities on Maker's side.  
In order to state Breaker's strategy precisely, we 
introduce a definition.  
We say that a color is {\em dead} if it is 
available on less than \( 6 (\log n)(\log_b np)\) vertices
on Maker's side $W_M$.
Breaker colors according to the following three simple rules:
\begin{enumerate}
\item Only color on Breaker's side $W_B$,
\item Do not use a dead color, and
\item If possible, respond to a move by Maker on Maker's side in 
kind (i.e. when Maker plays on Maker's side with a particular color then
Breaker's first choice is to play the same color on Breaker's side).
\end{enumerate}

We say that a color {\em escapes} if it is not dead
and Breaker stops playing this color
because he cannot choose a vertex on his side
that can be colored with this color. 
Note that it follows from the third rule for Breaker that the number of
times a color is played on Breaker's side is at least 
the number times it is played on
Maker's side as long as the color is neither dead nor has escaped.
Note further that there may be rounds when Breaker's move will
not be dictated by the rules above.  During these
rounds Breaker simply colors arbitrarily on Breaker's side.
We continue play until every color either dies or
escapes; that is, we play until Breaker cannot follow his 
coloring rules. Suppose that this happens after $\n_M \leq \n_B$ vertices have been
colored on Maker's and Breaker's sides respectively.
We will show that \whp\ Breaker will be in a winning position by this time.

Recall that $W_M,W_B$ denote Maker and Breakers' sides of the bipartition.  For \( X \subseteq 
W_B \) we let \( \bN(X) \) denote the set of vertices in \( W_M \) that have no neighbors in
\(X\) (i.e. \( \bN(X) = W_M \setminus N(X) \)).  Let
$$\l_0=\log_bn-\log_b\log n-\log_b\log_b np-\log_b3$$
and note that 
\begin{equation}\label{N3}
\l_0< \log_b np \hspace{1cm} \mbox{and} \hspace{1cm} n(1-p)^{\l_0}= 3 (\log n) (\log_b np).
\end{equation}

\begin{lemma}\label{lem9}
\Whp\ every subset $L\subseteq W_B$ 
of size $\ell \leq \l_0$ has at most
$2n(1-p)^{\ell }$ non-neighbors in $W_M$. 
\end{lemma}
\proofstart
Fix $L \subseteq W_B$ with $|L|=\ell$. The number of non-neighbors of $L$ in $W_M$ is 
distributed
as the binomial $B(n,(1-p)^\ell )$. Thus, by the Chernoff bounds,
$$\Pr(\exists L:|\bN(L)|\geq 2n(1-p)^{\ell})\leq\sum_{\ell =1}^{\l_0}
\binom{n}{\ell }e^{-n(1-p)^\ell /3}\leq
\sum_{\ell =1}^{\l_0}\bfrac{ne}{\ell }^\ell e^{-n(1-p)^\ell /3}.$$
Now if $\ell \leq \l_0<\log_b np$ then
$$\log(n^\ell e^{-n(1-p)^\ell /3})= \ell\log n-n(1-p)^\ell/3\leq \begin{cases}-\sqrt{n}& \text{ 
if }\ell \leq
  \sqrt{\log n}\\0&  \text{ if } \sqrt{\log n} \leq \ell \leq \l_0\end{cases}.$$

Therefore,
$$\Pr\big(\exists L:|\bN(L)|\geq
2n(1-p)^{\ell}\big)\leq\sum_{\ell =1}^{\sqrt{\log n}}\bfrac{e}{\ell }^\ell e^{-\sqrt{n}}+
\sum_{\ell =\sqrt{\log n}}^{\l_0}\bfrac{e}{\ell }^\ell =o(1).$$
\proofend

It follows from \eqref{N3}, Lemma~\ref{lem9} and the definition of a dead color that 
Breaker makes a rule based use of each color at most $\l_0$ times.
Using the fact that at least as many vertices will be colored on Breaker's side as on Maker's 
side and that they both had the same number of turns we conclude that the number of 
colored vertices at the point when Breaker stops satisfies
$$\n_M \leq \n_B \leq 2k\l_0 \leq \frac{n}{5\log n}.$$
Let \( c_1, \dots, c_t \) be the colors that escape.    
Let \( M_i,B_i \) be the sets of vertices with color \(c_i\) on Maker's and Breaker's sides respectively at the moment that
Breaker stops playing color \( c_i\) because he is forced to by the rules.  Let \( m_i = |M_i| \) and set
$$\alpha = \sum_{i=1}^t (1-p)^{m_i}. $$
Note first that
$$m_i\leq b_i=|B_i|\leq \l_0,\qquad i=1,2,\ldots,t.$$
Furthermore, because $b_i\leq \l_0$ we see that 
$$|\bN(B_i)|\leq 2n(1-p)^{b_i}\leq 2n(1-p)^{m_i}.$$
We consider two cases.

\noindent{\bf Case 1.}  \( \alpha < 1/6 \).

\noindent
The total number of vertices that can be colored on Maker's side is at most 
the sum of (i) the number of
vertices colored so far, (ii) the number of vertices that can be
colored with dead colors, and (iii) the number of vertices that can be
colored with escaped colors.  Hence the number of vertices that can be colored on
Maker's side is at most
$$\n_M + k\cdot 6 (\log n)(\log_bnp)+ 2n\sum_{i=1}^t(1-p)^{m_i} 
\leq o(n)+\frac{n}{10 (\log n)(\log_bnp)}\cdot 6(\log n)(\log_bnp)+\frac{n}{3}<n,$$
and therefore Maker can not complete the coloring of the graph.

\noindent{\bf Case 2.}   \( \alpha \geq 1/6 \).

\noindent
In this case \whp\ we arrive at a contradiction.  Let \( Z\) be the set of
\( \n_B \) vertices in $W_B$ that have been colored so far.  We have
\( |Z|=\n_B\leq n/(5\log n)\).  When a color $c_i$ escapes it is unavailable to 
vertices on Breaker's side. 
It follows that all vertices in $Y=W_B\setminus Z$ have at least one neighbor in
\( M_i \).

Let \( \cE \) be the event that we have such a configuration i.e. $t$
small sets, whose neighborhoods each covers almost all of $W_B$.  Fix the sets \(M_1, \dots, 
M_t\) and \(Y\).  Since $|Y|=(1-o(1))n$,
the probability that this collection of sets satisfies the condition is
\[ \left( \prod_{i=1}^t\left( 1 - (1-p)^{m_i} \right) \right)^{|Y|} 
\le \exp \left\{ - |Y| \sum_{i=1}^t (1-p)^{m_i} \right\}\leq e^{-n/7}. \]
The probability of the existence of any such configuration in our random model is at most
\[ \sum_{t=1}^k \brac{\sum_{\ell=0}^{\l_0}\binom{n}{\ell}}^t \binom{n}{n/(5\log n)}
e^{-n/7}\leq n^{\l_0k}e^{o(n)}e^{-n/7}\leq e^{n/10+ o(n) -n/7}=o(1) \]
\proofend

\subsection{Upper bound in Theorem \ref{th2}}
The proof here is essentially the same as for Theorem
\ref{th1}(a). Let the vertex bipartition be denoted $V_1,V_2$. Aside from modifying the statement of Lemma \ref{lem5}
to say that the sets in $\cC$ are contained in $V_i$ and $v\in
V_{3-i}$, the proof goes through basically unchanged.

\section{Concluding remarks and open problems}
In this paper we obtain upper and lower bounds on the game chromatic number of $G_{n,p}$
which differ only by a multiplicative constant. It would be very interesting to improve our 
result and  determine the asymptotic value of this parameter for random graphs.
Our results suggest that in fact the following should be true.

\begin{conj}
If $p\leq 1-\eta$ for some constant $\eta>0$ and $np \to\infty$,
then {\bf whp} 
$$\chi_g(G_{n,p})=(1+o(1))\frac{n}{\log_b np},$$
where $b=1/(1-p)$.
\end{conj}
We conjecture that
the game chromatic number of random bipartite graph $B_{n,p}$ has the same 
order of magnitude {\bf whp}.   We did not succeed in 
proving the correct lower bound.

As one final remark, there has been much work done on the concentration of the the chromatic 
number of $G_{n,p}$. None of this is applicable to $\chi_g(G_{n,p})$. Now $\chi_g(G_{n,p})$ 
should also be concentrated, but proving this may
require some new approaches to proving concentration.

\end{document}